\documentstyle[11pt]{article}
\input amssym.def

\newtheorem{thm}{Theorem}[section]
\newtheorem{lem}[thm]{Lemma}
\newtheorem{cor}[thm]{Corollary}
\newtheorem{pro}[thm]{Proposition}
\newtheorem{ex}[thm]{Example}

\newtheorem{defi}[thm]{Definition}

\newcommand{\gm }{\Gamma }
\newcommand{\lon }{\longrightarrow }
\newcommand{\be }{\begin{eqnarray*}}
\newcommand{\ee }{\end{eqnarray*}}

\newcommand{\per }{\backl }

\newcommand{\poidd }[2]{#1\gpd #2}

\def\nin{\noindent}

\newcommand{\pf}{\noindent{\bf Proof.}\ }
\newcommand{\qed}{\begin{flushright} $\Box$\ \ \ \ \ \
                  \end{flushright}}
\newcommand{\complex}{{\Bbb      C}}
\newcommand{\reals}{{\Bbb      R}}

\newcommand{\cinf}{C^{\infty}}

\newcommand{\id}{\mbox{\rm{id}}}
\newcommand{\backl}{\mathbin{\vrule width1.5ex height.4pt\vrule height1.5ex}}

\newcommand{\cald}{{\cal D}}

\newcommand{\calg}{{\cal G}}
\newcommand{\calh}{{\cal H}}

\newcommand{\calk}{{\cal K}}
\newcommand{\call}{{\cal L}}

\newcommand{\calu}{{\cal U}}
\newcommand{\calx}{{\cal X}}
\newcommand{\caly}{{\cal Y}}

\newcommand{\smalcirc}{\mbox{\,\tiny{$\circ $}\,}}     

\def\description label#1{\hfil\bf[#1]\hfil}
\newcommand{\sect}{\gm(A)}

\newcommand{\seccg}[1]{\gm(\wedge^{#1}AG)}                 

\newcommand{\anc}{a}
\newcommand{\ancd}{a_{*} }

\newcommand{\Ga}{\Gamma}
\newcommand{\pib}{\pi^{\#}}
\newcommand{\talp}{\widetilde{\alpha}}           
\newcommand{\tbet}{\widetilde{\beta}}            
\newcommand{\difft}{\frac{d}{dt}}

\def\LAgpd{${\cal LA}$-groupoid}

\def\Ad{\mathop{\rm Ad}}

\def\ST{\ \vert\ }

\def\tilalpha{\widetilde\alpha}
\def\tilbeta{\skew6\widetilde\beta}

\def\sdp{\mathbin{\hbox{$\mapstochar\kern-.3333em\times$}}}
\def\pds{\mathbin{\hbox{$\times\kern-.55em\mapstochar\,$}}}

\newcommand{\wed}{\mathbin{\lower1.5pt\hbox{$\scriptstyle{\wedge}$}}}

\let\Tilde=\widetilde
\let\Bar=\overline
\let\Ri=\overrightarrow
\let\Le=\overleftarrow
\def\ssri{{\stackrel{\rightarrow}{X}}}
\def\ssli{{\stackrel{\leftarrow}{X}}}
\def\ssriy{{\stackrel{\rightarrow}{Y}}}

\let\sol=\bullet
\def\chigh{{\raise1.5pt\hbox{$\chi$}}}
\let\phi=\varphi

\def\til0{\Tilde{0}}
\def\til1{\Tilde{1}}

\def\dminus{\raise2pt\hbox{\vrule height1pt width 2ex}\hskip3pt}

\def\pback#1{\mathbin{{{\lower1.2ex\hbox{$\times$}}\atop #1}}}

\def\ddt#1{\left.\frac{d}{dt}#1\right|_0}

\def\vlra{\hbox{$\,-\!\!\!-\!\!\!-\!\!\!-\!\!\!-\!\!\!
-\!\!\!-\!\!\!-\!\!\!-\!\!\!-\!\!\!\longrightarrow\,$}}

\def\gpd{\,\lower1pt\hbox{$\longrightarrow$}\hskip-.24in\raise2pt
             \hbox{$\longrightarrow$}\,}

\def\lgpd{\,\lower1pt\hbox{$\vlra$}\hskip-1.02in\raise2pt\hbox{$\vlra$}\,}

\def\llgpd{\,\lower1pt\hbox{$\vvlra$}\hskip-1.3in\raise2pt\hbox{$\vvlra$}\,}

\def\vgpd{\Bigg\downarrow\!\!\Bigg\downarrow}

\begin{document}

\title{{\bf Integration of Lie bialgebroids}
}

\author{KIRILL C. H. MACKENZIE\\
        School of Mathematics and Statistics\\
        University of Sheffield\\
        Sheffield, S3 7RH, England\\
        {\sf email: K.Mackenzie@sheffield.ac.uk}\\
        and \\
         PING XU \thanks {Research partially supported by NSF
        Grant DMS95-04913.}\\
Department of Mathematics\\
The Pennsylvania State University\\
University Park, PA 16802, USA\\
        {\sf email: ping@math.psu.edu }}

\date{{\sf December 21, 1997}}

\maketitle

\begin{abstract}
We prove that under certain mild assumptions a Lie bialgebroid integrates
to a Poisson groupoid. This includes, in particular, a new proof of the
existence of local symplectic groupoids for any Poisson manifold, a
theorem of Karasev and of Weinstein.
\end{abstract}

\section{Introduction}
A symplectic realization of a Poisson manifold $P$ is a Poisson map from
a symplectic manifold $X$ to $P$ which is a surjective submersion.
The idea of finding symplectic realizations for degenerate Poisson brackets
can be traced back to Lie, who used the name  ``function group" \cite{Lie}.
Lie proved that such a realization always exists locally for any Poisson
manifold of constant rank. The local existence theorem of symplectic
realizations for general Poisson manifolds was proved by Weinstein in 1983
\cite{Weinstein:1983}. The proof was highly nontrivial and used the local
structure theorem for Poisson manifolds. Subsequently, in 1987, Karasev
\cite{Karasev:1987} and Weinstein \cite{Weinstein:1987} proved independently
the existence of a global symplectic realization for any Poisson manifold.
In fact, they found that by a suitable choice, such a realization admits
{\bf automatically} a local groupoid structure which is compatible with the
symplectic structure in a certain sense. The global form of this notion is
what is now called a symplectic groupoid. Symplectic groupoids have their
own origin in quantization theory \cite{KarasevM}. However, it has been
quite mysterious why the groupoid structure and symplectic structure enter
into the picture of a Poisson manifold in such a compatible and striking
manner.

On the other hand, Poisson groups have been intensively studied as a
classical limit of quantum groups. The theory of Poisson groups established
a precise relation between Poisson structures on the groups and their
infinitesimal invariants, Lie bialgebras. In order to understand symplectic
groupoids using the techniques of Poisson group theory and to unify both
theories in a general framework, Weinstein in 1988 introduced the notion
of Poisson groupoid \cite{Weinstein:1988}. Lie bialgebroids were introduced
and studied by the present authors \cite{MackenzieX:1994} in 1994 as the
infinitesimal invariants
of Poisson groupoids: given a Poisson groupoid $G$, the Lie algebroid of the
underlying Lie groupoid, together with the Lie algebroid structure on the
dual $A^*G$ \cite{Weinstein:1988}, form a Lie bialgebroid. Lie bialgebroids
are found to be connected with various subjects in Poisson geometry ranging
from Poisson-Nijenhuis structures to Dirac structures; see, for example,
\cite{Kosmann-Schwarzbach:1995}, \cite{Kosmann-Schwarzbach:1996},
\cite{LiuWX:1997}, \cite{LiuWX:Dirac}, \cite{Mackenzie:SQ2}. However, it has
remained an unsettled problem whether an arbitrary Lie bialgebroid can be
integrated to a Poisson groupoid.

In this paper, we give an affirmative answer to this question. We prove that
a Lie bialgebroid structure on the Lie algebroid of a (suitably
simply-connected) Lie groupoid can be integrated to give a Poisson groupoid
structure on the underlying groupoid. This result extends the well-known
result that a Lie bialgebra (of finite dimension over $\reals$ or $\complex$)
is the Lie bialgebra of a Poisson group \cite{Drinfeld:1983}, \cite{LuW:1990}.
At the other extreme, it also shows that if a Poisson manifold $P$ has a
cotangent Lie algebroid which integrates to a Lie groupoid $G\gpd P$, then
$G$ has a symplectic groupoid structure integrating the Poisson structure of
$P$. This is a large part of the local integrability of Poisson manifolds
\cite{Karasev:1989}, \cite{Weinstein:1983}, \cite{Weinstein:1987}.
In particular, we obtain as a consequence a new proof of the existence of
local symplectic groupoids for general Poisson manifolds.

Within this general framework, the geometric origin of the symplectic and
groupoid structures on a symplectic groupoid becomes transparent. Given a
Poisson manifold $P$, its cotangent bundle $T^* P$ carries a Lie algebroid
structure (see Section~\ref{sect:sg} for the precise definition). Now
if this integrates to a Lie groupoid $\Gamma$ (assumed $\alpha$--simply
connected), the canonical Lie algebroid structure on the dual of $T^*P$,
namely the tangent bundle $TP$, induces a Poisson structure on $\Gamma$,
which in this case is symplectic. The compatibility condition
between the two Lie algebroid structures then ensures the compatibility
condition between the groupoid and symplectic structures, and so $\Gamma$ is
a symplectic groupoid.

The integrability of Lie bialgebras reduces essentially to the lifting, or
integration, of  Lie algebra 1-cocycles, but for general Lie bialgebroids
this approach is not available: there is no satisfactory adjoint
representation for a general Lie algebroid, and if one were to treat the Lie
algebroids as infinite dimensional Lie algebras, then the problems would be
at least as great as those resolved here.

In Theorem 6.2 of \cite{MackenzieX:1994}, we proved that a Lie algebroid
$A$, whose dual $A^*$ also has a Lie algebroid structure, is a Lie
bialgebroid if and only if a certain map $\Pi\colon T^*A^*\lon TA$ is a Lie
algebroid morphism. Here $T^*A^*$ is the cotangent Lie algebroid for
the Poisson structure on $A^*$ induced by $A$, and $TA$ is the tangent
Lie algebroid structure on base $TP$, of $A\lon P$ \cite[\S5]{MackenzieX:1994}.
This result now allows us, after some work, to reduce the integrability problem
to the integration of Lie algebroid morphisms. This approach is not, in fact,
so very far from the integration of cocycles---the standard proof of the
integrability of cocycles proceeds by treating them as morphisms into
semi-direct products. Here, however, none of the Lie algebroid structures need
be a semi-direct product.

Unlike the case of Lie algebras, a Lie algebroid need not arise from a Lie
groupoid. If $G$ is a Poisson groupoid and the Lie algebroid dual $A^*G$
integrates to an $\alpha$-simply connected groupoid $G^*$, then the results
of the present paper make $G^*$ into a Poisson groupoid also. The further
relations between a Poisson groupoid and its dual will be investigated
elsewhere.

The results of the present paper settle, we believe, any remaining doubt
that the concept of Lie bialgebroid is the correct infinitesimal form of
the concept of Poisson groupoid. This is an important point, in view of the
complexity of the work on notions of double for Lie bialgebroids
\cite{LiuWX:1997}, \cite{Mackenzie:SQ2}.

We begin in \S2 by giving some preliminary material concerning affine
multivector fields on Lie groupoids. In \S3 we recall the basic definitions
and main results from \cite{MackenzieX:1994}. The main integrability theorem
is proved in \S4. In \S5 we consider the case of symplectic groupoids proving,
in particular, that if the cotangent Lie algebroid of a Poisson manifold
integrates to an $\alpha$-simply connected groupoid, then the groupoid has
a natural structure of symplectic groupoid. Thus the existence of a local
symplectic groupoid for a Poisson manifold follows as a consequence.
Finally in \S\ref{sect:app} we give a full proof of the integrability result
for Lie algebroid morphisms on which the main results depend.

We have minimized the repetition of material from \cite{MackenzieX:1994},
and so have used the same notation conventions.

We are once again very grateful to Alan Weinstein for conversations over an
extended period on the material of this paper. We also thank the Isaac Newton
Institute at Cambridge and the organizers of the workshop on Symplectic
Geometry. The second author wishes to thank IHES and Max-Planck-Institut
for their hospitality while part of the work was being done.

\section{Affine multivector fields}

Throughout this section, we fix a Lie groupoid $\poidd{G}{M}$ and denote its
Lie algebroid by $AG$. We follow the conventions of \cite{MackenzieX:1994};
in particular, for $g,h\in G$ the product $gh$ is defined if
$\alpha g = \beta h$.

We recall the exponential map for a Lie groupoid \cite{Mackenzie:LGLADG},
\cite{CDW}. Given $X\in\gm AG$, the flows $\phi_t$ of the corresponding
right invariant vector field $\Ri{X}$ are left translations. Assume for
convenience that $\phi_t$ is global and define $\exp{tX}\colon M\lon G$
by $\exp{tX(m)} = \phi_t(1_m)$. Then $\exp{tX}$ is a section of
$\alpha\colon G\lon M$, and $\beta\smalcirc\exp{tX} \colon M\lon M$ is the
flow of $a(X)\in\calx(M)$. Call any section $\calk \colon M\lon G$
of $\alpha$ for which $\beta\smalcirc\calk$ is a diffeomorphism a {\em bisection}
of $G$ ({\em admissible section} in \cite{Mackenzie:LGLADG}). Then
$L_{\calk}(g) = \calk(\beta g)g$ is the left translation corresponding to
$\calk$, and $R_\calk(g) = g\calk((\beta\smalcirc\calk)^{-1}(\alpha g))$ is the
right translation. We often denote the tangents of $L_\calk$ and $R_\calk$
by the same symbols. By $Ad_\calk$ we denote the groupoid automorphism
$L_\calk\smalcirc R^{-1}_\calk$. It is clear that $Ad_\calk$ leaves $M$
invariant, and its restriction to $M$ is the map $\beta\smalcirc\calk$.
The set of all (global) bisections $\calg(G)$ forms a group under
$\calk_1\smalcirc\calk_2(m) = \calk_1( \Ad_{\calk_{2}}(m)) \calk_2(m)$.
In general $\exp{tX}$ is a local bisection (in an evident
sense), and is only defined for small $t$.

Alternatively, one may identify a bisection with its image, in which case a
bisection is a submanifold of $G$ for which the restrictions of both $\alpha$
and $\beta$ are diffeomorphisms \cite{CDW}. Then $\exp{tX}$, for $X\in\gm AG$,
is the submanifold of $G$ obtained by flowing the identity space $M$ under
the flow $\phi_t$ of $\Ri{X}$. We will use both points of view in what
follows. The following formulas are frequently used in the paper.
$$
\alpha (\exp{tX}(m)) =m, \qquad
\beta (\exp{tX}(m)) = Ad_{\exp{tX}}(m).
$$
\begin{defi}
A multivector field $D$ on $G$ is {\em affine} if for any $x, y\in G$ such
that $\alpha(x)= \beta(y)=m$ and any bisections $\calx, \caly$ through the
points $x, y$, we have
\begin{equation}
\label{eq:multi}
D(xy) = R_{\caly }D(x) + L_{\calx } D(y) - R_{\caly }L_{\calx }D (1_m).
\end{equation}
\end{defi}

To appreciate this definition, recall that $TG$ inherits a groupoid
structure on base $TM$ from $G\gpd M$ with source map $T(\alpha)$,
target map $T(\beta)$, and composition $X\sol Y = T(\kappa)(X,Y)$, where
$\kappa$ is the composition in $G$. One of us proved \cite[2.4]{Xu:1995}
that if $X\in T_x(G)$ and $Y\in T_y(G)$ have
$T(\alpha)(X) = T(\beta)(Y)  = W$, then this composition may also be given by
\begin{equation}                             \label{eq:Xumult}
X\sol Y =
T(L_{{\cal X}})(Y) + T(R_{{\cal Y}})(X) -
            T(L_{{\cal X}})T(R_{{\cal Y}})(T(1)(W))
\end{equation}
where ${\cal X}, {\cal Y}$ are any (local) bisections of $G$ with
${\cal X}(\alpha x) = x$ and ${\cal Y}(\alpha y) = y$.

For the case of affine multivector fields on groups, see
\cite[\S4]{Weinstein:1990}. Affine multivector fields arise in many natural
ways. For instance, for any $K\in \gm (\wedge^{k}AG)$, the multivector
field $\Ri{{K}}-\Le{K}$, the difference between the right and left
translations of $K$, is an affine multivector field on $G$. We deal
elsewhere \cite{MackenzieX:1998} with the special features of affine
vector fields.

The following theorem gives a very useful characterization of affine
multivector fields. See \cite{LuW:1990} and \cite{Weinstein:1990}
for the case of groups.

\begin{thm}                                     \label{thm:affine}
Let $G$ be $\alpha$-connected. For a multivector field $D$ on $G$, the
following statements are equivalent:
\begin{enumerate}
\item $D$ is affine;
\item For any right (left) invariant vector field $\Ri{X}$ (respectively
$\Le{X}$), the Lie derivative $L_{{\ssri}} D$ ($L_\ssli D$) is right (left)
invariant;
\item for any left invariant vector field $\Le{X}$ and right invariant vector
field $\Ri{Y}$, $L_\ssriy L_\ssli D = 0$.
\end{enumerate}
\end{thm}

\pf
The equivalence between (ii) and (iii) is quite evident.

$\mbox{(i)}\Longleftrightarrow \mbox{(ii)}$
Our proof here follows that of Theorem 3.1 in \cite{Xu:1995}.

Suppose that $D$ is an affine multivector field. It suffices to show that
statement (ii) holds for any compactly supported $X\in \gm_{c} (AG)$,
since the evaluation of $L_\ssri D$ at any particular point only depends
on the local germ of $X$. Therefore, without loss of generality, we can
assume that $\Ri{X}$ is a complete vector field on $G$.
Let $\calx_{t}=\exp {tX}\in \calg (G)$, the one parameter family
of bisections generated by $X\in \gm_{c}(AG)$. Suppose that $\caly$ is
any bisection through a point $y\in G$ with $\beta  (y)=m$. Let
$x_{t} = \calx_t(m)$, the flow generated by $\Ri{X}$ with initial point
$1_m$. According to Equation (\ref{eq:multi}), we have
$$
D(x_{t} y) = R_{\caly } D (x_{t}) + L_{\calx_{t} }D (y) -
                                 L_{\calx_{t} }R_{ \caly }D (1_m),
$$
which is equivalent to
\begin{equation}
\label{eq:family}
L_{ \calx_{t}^{-1} } D (x_{t} y) =  (R_{\caly } \smalcirc L_{\calx_{t}^{-1} })
D (x_{t})+ D (y) -R_{\caly } D (1_m).
\end{equation}
Here, both sides of Equation (\ref{eq:family}) are elements in
$\wedge^{k}T_{x}G$, where $k$ is the degree of $D$.
Taking the derivative at $t=0$, one obtains that
\begin{equation}
\label{eq:lie-derivative}
(L_\ssri D )(y) = R_{\caly } [(L_\ssri  D)(1_m)].
\end{equation}
This implies that $L_\ssri D$ is right-invariant, according to Lemma~3.2 in
\cite{Xu:1995}.

Conversely, integrating Equation (\ref{eq:lie-derivative}),
one gets immediately that
\begin{equation}
\label{eq:inter}
L_{\calx_{t}^{-1} } D (x_{t} y)-D (y) =
         R_{\caly}[ L_{  \calx_{t}^{-1} } D (x_{t}) -D (1_m)].
\end{equation}
Thus Equation (\ref{eq:multi}) follows if $\caly = \exp{X}$ for
any $X\in \gm_{c}(AG)$. In the general case it follows by applying Equation
(\ref{eq:inter}) repeatedly and the following result.

\begin{lem}                                          \label{lem:enough}
For any $x\in G$ and any bisection $\calk$, there exist
$X_{1}, \ldots , X_{n}\in \gm_{c}(AG)$
such that $\tilde{\calk}=\exp{X_{1}}\ldots \exp{X_{n}}$ has the value $x$
at $\alpha x$ and $\calk_{*}|_{\alpha x}=\tilde{\calk}_{*}|_{\alpha x}$.
\end{lem}

\pf
Let ${\cal U}$ denote the set of values $\exp{tX}(m)$ as $X$ ranges through
$\gm_c(AG)$, $t$ ranges through $\reals$, and $m$ through $M$.
Clearly, the statement holds if $x\in \calu$. It is easy
to verify that the intersections of ${\cal U}$ with each $\alpha$-fibre
is open and it follows by a modification of a standard argument
\cite[II\S3]{Mackenzie:LGLADG} that the subgroupoid $H$ generated by
${\cal U}$ has each $\alpha$-fibre open in the corresponding
$\alpha$-fibre of $G$. Since $G$ is $\alpha$-connected it follows that
$H = G$.
\qed

Any bivector field $D$ induces a map $D^{\#}: T^* G\lon TG$ defined by
$\langle\omega_2,D^\#(\omega_1)\rangle =
\langle D,\omega_1\wed\omega_2\rangle$. It is well known that both
$T^* G$ and $TG$ have Lie groupoid structures \cite{CDW},
\cite{Pradines:1988}, that for $TG$ being given above (\ref{eq:Xumult}).
For $T^*G$ we use the conventions of \cite[\S7]{MackenzieX:1994}:
we take the source $\tilalpha$ and target $\tilbeta$ to be given by
\begin{equation}                    \label{eq:T*G}
\Tilde{\alpha}(\omega)(X) = \omega(T(L_g)(X - T(1)(a(X)))),\qquad
\Tilde{\beta}(\omega)(Y) = \omega(T(R_g)(Y)),
\end{equation}
where $\omega\in T^*_gG,\ X\in A_{\alpha g}G$ and $Y\in A_{\beta g}G$. If
$\theta\in T^*_hG$ and $\Tilde{\alpha}(\theta) = \Tilde{\beta}(\omega)$
then $\alpha h = \beta g$ and we define $\theta\sol\omega\in T^*_{hg}G$ by
$$
(\theta\sol\omega)(Y\sol X) = \theta(Y) + \omega(X),
$$
where $Y\in T_hG,\ X\in T_gG$. The identity element
$\Tilde{1}_\phi \in T^*_{1_m}G$ corresponding to $\phi\in A^*_mG$ is defined
by $\Tilde{1}_\phi(T(1)(x) + X) = \phi(X)$ for $X\in A_mG,\ x\in T_m(M)$.

We can now give the following criterion for $D$ to be affine, which will be
important in the proof of Theorem \ref{thm:main}. An immediate consequence of
this criterion is that the Poisson tensor on a Poisson groupoid is affine,
one of the main results proved in \cite{Xu:1995}. In fact, our proof
here is essentially borrowed from that in \cite{Xu:1995}.

\begin{pro}                                \label{pro:D}
Let $D$ be a bivector field on $G$.
If
\begin{equation}
\matrix{&&D^\#&&\cr
	&T^*G&\vlra&TG&\cr
		&&&&\cr
			&\vgpd&&\vgpd&\cr
				&&&&\cr
               &A^*G&\vlra&TP&\cr
						&&a_*&&\cr}
\end{equation}
is a Lie groupoid morphism, then the bivector field $D$ is affine.
\end{pro}

\pf By $\Lambda \subset G\times G\times G$, we denote the
graph of groupoid multiplication, and by $\Omega$, we denote
the subset of $G \times G \times G  \times G$
consisting of all elements $(z, y, x, w)$ such that $w=yz^{-1}x$.
Weinstein  \cite{Weinstein:1990} calls $\Omega$  the {\em affinoid diagram}
corresponding to the groupoid $G$.  The graph of multiplication of the
groupoid $T^{*}G\gpd A^*G$ is $\bar{N}^{*}\Lambda $, which is the
subset of $ T^{*}(G\times G \times G)$ obtained from the conormal bundle
$N^{*}\Lambda$ by multiplying the cotangent vectors in the last factor
by $-1$. Thus, by assumption, $D^{\#} \bar{N}^* \Lambda \subset T\Lambda$.
Following the proof of Theorem 4.5 in \cite{Weinstein:1990},
we can show that $D^{\#} N^*\Omega  \subset T\Omega $.
For any  $x, y\in G $ such that $\alpha (x)=\beta (y)=m$, it is clear that
$(z,y,x,1_m)$, with  $z=xy$, is an element of $\Omega$. For any
$\xi \in T_{z}G $, it follows from Lemma 2.6 in \cite{Xu:1995} that
$(-\xi, L^{*}_{\calx}\xi, R^{*}_{\caly}\xi, - L^{*}_{\calx} R^{*}_{\caly}\xi)$
is a conormal vector to $\Omega$. Therefore, for any $\xi, \eta \in T_{z}G$,
we have
$$
-D (z) (\xi, \eta) +D  (y)( L^{*}_{\calx}\xi, L^{*}_{\calx}\eta )
+D (x) (  R^{*}_{\caly}\xi ,  R^{*}_{\caly}\eta )
-D (w) ( L^{*}_{\calx} R^{*}_{\caly}\xi ,  L^{*}_{\calx} R^{*}_{\caly}\eta )
=0.
$$
This implies Equation (\ref{eq:multi}) immediately.
\qed

A bivector field $D$ for which $D^\#$ is a morphism might be called
{\em multiplicative}. Not all affine bivector fields are multiplicative.

\begin{pro}
The Schouten bracket of affine multivector fields is still affine.
\end{pro}

\pf Suppose that $D_{1}$ and $D_{2}$ are affine multivector fields on $G$.
For any $X, Y\in \gm (AG)$, we have
$$
L_\ssriy L_\ssli  [D_{1}, D_{2}]=[L_\ssriy L_\ssli  D_{1}, D_{2}]
+[L_\ssli  D_{1}, L_\ssriy  D_{2}]+[L_\ssriy D_{1}, L_\ssli   D_{2}]
+[D_{1}, L_\ssriy L_\ssli  D_{2}].
$$
Each summand on the right hand side is easily seen to be zero
according to Theorem \ref{thm:affine}, so it follows that
$[D_{1}, D_{2}]$ is affine.
\qed

As in the case of groups, the derivative of an affine $k$-vector field
$D$ can be introduced, and is a map $dD : \gm (AG)\lon \gm (\wedge^k AG)$
defined as follows. For any $X\in \gm (AG)$, $dD (X)$ is defined to be the
element in $\gm (\wedge^k AG)$ whose right translation is $L_\ssri D$. It
is easy to see that if $\gm(AG)$ is considered as an infinite dimensional
Lie algebra with $\gm (\wedge^k AG)$ considered as a $\gm (AG)$-module in
a natural way, then $dD$ may be considered a Lie algebra 1-cocycle.
It is not clear in general whether such a Lie algebra 1-cocycle can be
lifted to an affine multivector field. However, the following
theorem indicates that if it exists, then it is unique.

\begin{thm}
\label{thm:dD0}
Let $D$ be an affine multivector field on an $\alpha$-connected Lie
groupoid $G$. Then $D$ is zero if and only if $D$ vanishes on the unit
space $M$ and $dD=0$.
\end{thm}

\pf
Let $X\in \gm_{c} (AG)$ be any compactly supported section and let
$\calx_{t}=\exp {tX}\in \calg (G )$ be the one parameter family of bisections
generated by $X$. Fixing any $m\in M$, let $x_{t}=\exp {tX} (m)$, and let
$f_{m}(t)$ be the vector in $\wedge^{k}T_{m}G$ given by
$f_{m}(t)=L_{\calx_{t}^{-1}}D(x_{t})$, where $k$ is the degree of $D$.
For any $t, s\in \reals$,
\be
f_{m}(t+s)&=&L_{\calx_{t+s}^{-1}}D(x_{t+s})\\
&=& L_{\calx_{t+s}^{-1}}D[\exp {sX}(v)\cdot \exp {tX}(m )] \\
&=& L_{\calx_{t+s}^{-1}}[R_{\calx_{t}}D(\exp {sX}(v))+L_{\calx_{s}}
         D(\exp {tX}(m ))]\\
&=& Ad_{\calx_{t}^{-1}} L_{\calx_{s}^{-1}}[D(\exp {sX}(v ))]+
L_{\calx_{t}^{-1}}D(\exp {tX} (m)) \\
&=& Ad_{\calx_{t}^{-1} } f_{v}(s)+ f_{m}(t),
\ee
where $v=\beta [\exp {tX}(m)]=Ad_{\calx_{t}}m$, and
in the third equality  we have used Equation (\ref{eq:multi}).

By taking the derivative with respect to $s$ at $0$, it follows immediately
that
$$
\difft f_{m}(t)=(Ad_{\calx_{t}}^{-1} )_{*}(L_\ssri  D)(v)=0.
$$
Therefore, $f_{m}(t)=0$ for all $t$ since $ f_{m}(0)=0$. This shows that
$D(\exp {tX}(m))=0$ for all $t$. Since any element in $G$ can be written
as a product of elements of the form $\exp{X}(m)$, it thus follows that
$D$ is identically zero on $G$, again by Equation (\ref{eq:multi}).
\qed

\section{Poisson groupoids and Lie bialgebroids}

In this section we briefly review some material from \cite{MackenzieX:1994}.

The concept of Lie bialgebroid can be defined in terms of a cocycle-type
condition, using a generalized Schouten calculus, but in this paper we
will be mainly concerned with the following equivalent characterization
in terms of morphisms of Lie algebroid structures.

\begin{defi}[\protect{\cite[3.1]{MackenzieX:1994}}]
\label{defi:bialg}
Suppose that $A\lon P$ is a Lie algebroid, and that its dual bundle
$A^{*}\lon P$ also carries a Lie algebroid structure. Then $(A, A^{*})$ is
a {\em Lie bialgebroid} if for any $X,\ Y\in \sect$,
\begin{equation}
d_*[X, Y] = L_{X}d_*Y-L_{Y}d_*X.
\label{eq:bialg}
\end{equation}
\end{defi}

For an alternative treatment of this definition, see Kosmann--Schwarzbach
\cite{Kosmann-Schwarzbach:1995}.

\begin{thm}[\protect{\cite[6.2]{MackenzieX:1994}}]
\label{thm:bialgebroid-equ}
Suppose that $q\colon A\lon P$ is a Lie algebroid such that its dual vector
bundle $q_*\colon A^{*}\lon P$ also has a Lie algebroid structure. Let
$\anc, \ancd$ be their anchors. Then $(A, A^{*})$ is a Lie bialgebroid if
and only if
\begin{equation}                         \label{eq:bialgebroid-equ}
\matrix{&\Pi &\cr
   T^*(A^*)&\vlra&TA\cr
      &&\cr
         \Bigg\downarrow&&\Bigg\downarrow\cr
            &&\cr
               A^*&\vlra&TP\cr
                  &a_*&\cr}
                  \end{equation}
is a Lie algebroid morphism, where the domain $T^{*}(A^{*})\lon A^{*}$ is
the cotangent Lie algebroid induced by the Poisson structure on $A^{*}$,
the target $TA\lon TP$ is the tangent prolongation of $A$, and
$\Pi \colon T^{*}(A^{*})\lon TA$ is the composition of the isomorphism
$R\colon T^{*}A^{*} \lon T^{*}A$ described below with
$\pi_{A}^{\#}: T^{*}A\lon TA$.
\end{thm}

We recall the structures used in this theorem. Given any vector bundle
$q\colon A\lon P$, the map $T(q)\colon TA\lon TP$ has a vector bundle
structure obtained by applying the tangent functor to the operations in
$A\lon P$. The operations in $TA\lon TP$ are consequently vector bundle
morphisms with respect to the tangent bundle structures in $TA\lon A$
and $TP\lon P$ and so $TA$ with these two structures is a double vector
bundle which we call the {\em tangent double vector bundle of} $A\lon P$
(see \cite[\S1]{Mackenzie:1992} and references given there). If
$q\colon A\lon P$ is a Lie algebroid then there is a Lie algebroid
structure on $T(q)\colon TA\lon TP$ defined in \cite[5.1]{MackenzieX:1994}
with respect to which $p_A\colon TA\lon A$ is a Lie algebroid morphism
over $p_P\colon TP\lon P$; we now call this the {\em tangent prolongation}
of $A\lon P$.

For a general vector bundle $q\colon A\lon P$, there is also a double vector
bundle
\begin{equation}                         \label{diag:T*A}
\matrix{&&r_A&&\cr
        &T^*A&\vlra&A^*&\cr
        &&&&\cr
   c_A&\Bigg\downarrow&&\Bigg\downarrow&q_*\cr
        &&&&\cr
        &A&\vlra&P&\cr
        &&q&&\cr}
\end{equation}
where $c_A$ is the usual cotangent bundle. For the structure on $r_A$ see
\cite[p.430]{MackenzieX:1994}. Elements of $T^*A$ can be represented locally
as $(\omega,X,\phi)$ where $\omega\in T^*_uP,\ X\in A_u,\ \phi\in A^*_u$ for
some $u\in P$. In these terms a canonical map $R\colon T^*A^*\lon T^*A$ can be
defined by $R(\omega,\phi,X) = (-\omega,X,\phi)$; for an intrinsic definition
see \cite[5.5]{MackenzieX:1994}. This $R$ is an isomorphism of double vector
bundles preserving the side bundles; that is to say, it is a vector bundle
morphism over both $A$ and $A^*$.

If $A^*$ has a Lie algebroid structure then the dual Poisson structure $\pi_A$
on $A$ has associated map $\pi_A^\#\colon T^*A\lon TA$ which is a morphism of
double vector bundles over $a_*\colon A^*\lon TP$ and $\id_A$. In a Lie
bialgebroid $(A,A^*)$ the same is consequently true of $\Pi$.

A {\em Poisson groupoid} (Weinstein \cite{Weinstein:1988}) is a Lie groupoid
$G\gpd P$ together with a Poisson structure $\pi_G$ on $G$ such that the
graph of the groupoid multiplication
$\Lambda = \{(h,g,hg)\ST\alpha h = \beta g\}$ is a
coisotropic submanifold of $G\times G\times\Bar{G}$. Any Poisson manifold
$P$ gives rise to a Poisson groupoid $\Bar{P}\times P\gpd P$ where $\Bar{P}$
is $P$ with the opposite structure and the groupoid structure is
$(w,v)(v,u) = (w,u)$. Any Poisson Lie group is of course a Poisson
groupoid. We consider examples further below.

It was shown in \cite{Weinstein:1988} that the manifold of identity elements
of a Poisson groupoid $G$ is coisotropic in $G$, and its conormal bundle
$N^*(P)$ thereby acquires a Lie algebroid structure. This conormal bundle
may be identified with $A^*G$, the dual vector bundle of $AG$, in a standard
way, and we will always take $A^*G$ with this Lie algebroid structure.
Denote the anchor of $A^*G$ by $a_*$.

In this paper we will use the following equivalent condition
\cite{AlbertD:1991ssr}, \cite[8.1]{MackenzieX:1994} repeatedly.

\begin{pro}
Let $G\gpd P$ be a Lie groupoid with a Poisson structure $\pi_G$. Then $G$
is a Poisson groupoid with respect to $\pi_G$ if and only if
$\pi_G^\#\colon T^*G\to TG$ is a morphism of Lie groupoids over some map
$a_*\colon A^*G\to TP$ (which is then the anchor of the Lie algebroid
dual).
\end{pro}

In a Poisson groupoid $G\gpd P$ we can therefore apply the Lie functor to
$\pi_G^\#$ and obtain a morphism of Lie algebroids
$A(\pi_G^\#)\colon AT^*G\to ATG$.

For any Lie groupoid $G\gpd M$, the tangent bundle projection
$p_G\colon TG\to G$ is a groupoid morphism over $p_M\colon TM\to M$ and
applying the Lie functor gives a canonical morphism $A(p_G)\colon ATG\to AG$.
This acquires a vector bundle structure by applying $A$ to the operations in
$TG\lon G$. This yields a system of vector bundles
\begin{equation}                         \label{diag:ATG}
\matrix{&&q_{TG}&&\cr
        &ATG&\vlra&TM&\cr
        &&&&\cr
  A(p_G)&\Bigg\downarrow&&\Bigg\downarrow&p_M\cr
        &&&&\cr
        &AG&\vlra&M&\cr
        &&q_G&&\cr}
\end{equation}
in which $ATG$ has two vector bundle structures, the maps defining each
being morphisms with respect to the other; that is to say, $ATG$ is a
double vector bundle.

Associated with the vector bundle $q_G\colon AG\lon M$ is the tangent double
vector bundle
\begin{equation}                         \label{diag:TAG}
\matrix{&&T(q_G)&&\cr
        &TAG&\vlra&TM&\cr
        &&&&\cr
  p_{AG}&\Bigg\downarrow&&\Bigg\downarrow&p_M\cr
        &&&&\cr
        &AG&\vlra&M.&\cr
        &&q_G&&\cr}
\end{equation}
It is shown in \cite[7.1]{MackenzieX:1994} that there is a canonical map
$$
j_G\colon TAG\to ATG
$$
which is an isomorphism of double vector bundles
preserving the side bundles. This $j_G$ is a restriction of the canonical
involution on $T^2G$.

Similarly, the cotangent groupoid structure $T^*G\gpd A^*G$ is defined by
maps which are vector bundle morphisms and, reciprocally, the operations in
the vector bundle $c_G\colon T^*G\lon G$ are groupoid morphisms. Taking the
Lie algebroid of $T^*G\gpd A^*G$ we get a double vector bundle
\begin{equation}                          \label{diag:AT*G}
\matrix{&&q_{T^*G}&&\cr
        &AT^*G&\vlra&A^*G&\cr
        &&&&\cr
  A(c_G)&\Bigg\downarrow&&\Bigg\downarrow&q_*\cr
        &&&&\cr
        &AG&\vlra&M,&\cr
        &&q_G&&\cr}
\end{equation}
where the vector bundle operations in $AT^*G\lon AG$ are obtained by applying
the Lie functor to those in $T^*G\lon G$.

It follows from the definitions of the operations in $T^*G\gpd A^*G$ that
the canonical pairing
$F\colon T^*G *_{G} TG\lon\reals,\ (\omega,X)\mapsto\langle\omega,X\rangle,$
can be considered a groupoid morphism into the the additive group(oid)
$\reals$. Here
$$
T^*G *_{G} TG =
\{(\omega,X)\ST \omega \in T^*G,\ X\in TG \mbox{ such that }
c_G(\omega) = p_G(X)\}
$$
is the pullback groupoid of $c_G$ and $p_G$; it has
base $A^*G *_M TM$. Hence $F$ induces a Lie algebroid morphism
$A(F)\colon AT^*G *_{AG} ATG\lon\reals$, where $AT^*G *_{AG} ATG$ is the
pullback Lie algebroid \cite[\S1]{HigginsM:1990a}. As noted in
\cite[7.2]{MackenzieX:1994}, $A(F)$ is nondegenerate, and so induces an
isomorphism of double vector bundles $i_G\colon AT^*G\lon A^\sol TG$, where
$A^\sol TG$ is the dual of $ATG\lon AG$. Now dualizing
$j_G\colon TAG\lon ATG$ over $AG$, we define
$$
j_G' = j_G^*\smalcirc i_G\colon AT^*G\lon T^*AG;
$$
this is an isomorphism of double vector bundles preserving the side bundles.
The following result is now immediate.

\begin{pro}
\label{pro:f-cocycle}
For $\nu\in AT^*G$ and $\xi\in ATG$ with $A(c_G)(\nu) = A(p_G)(\xi)$
we have
$$
A(F)(\nu,\xi) = \langle j_G'(\nu),j_G^{-1}(\xi)\rangle,
$$
where the pairing on the right is the canonical pairing between $T^*AG$ and
$TAG$.
\end{pro}

It is proved in \cite[7.3]{MackenzieX:1994} that the composition
$(j_G')^{-1}\smalcirc R$ is equal to the canonical isomorphism
$s\colon T^*A^*G\lon AT^*G$ arising from the symplectic groupoid structure
on $T^*G$, that is,
\begin{equation}
\label{eq:s}
s=(j_G')^{-1}\smalcirc R.
\end{equation}
In particular, $s$ is an isomorphism of double vector bundles from
\begin{equation}                          \label{diag:T*A*G}
\matrix{&&c_{A^*G}&&\cr
        &T^*A^*G&\vlra&A^*G&\cr
        &&&&\cr
r_{A^*G}&\Bigg\downarrow&&\Bigg\downarrow&q_*\cr
        &&&&\cr
        &AG&\vlra&M&\cr
        &&q_G&&\cr}
\end{equation}
to (\ref{diag:AT*G}) which preserves the side bundles.

Returning to the Poisson groupoid $G\gpd P$, we now have the morphism
$$
j_G^{-1}\smalcirc A(\pi_G^\#)\smalcirc {j'_G}^{-1}\colon T^*AG\to TAG
$$
and it is proved in \cite[\S8]{MackenzieX:1994} that this is equal to
$\pi_{AG}^\#$, the Poisson tensor map for the Poisson structure on $AG$
dual to the Lie algebroid structure on $A^*G$. It follows that
$\pi^\#_{AG}\smalcirc R = j_G^{-1}\smalcirc A(\pi_G^\#)\smalcirc s$ is a morphism of
Lie algebroids over $a_*$, and hence that $(AG, A^*G)$ is a Lie bialgebroid.
(An alternative proof of this result is given in \cite[3.5]{Xu:1995}.)
The integrability proof of \S\ref{sect:ilb} will consist essentially
of reversing these steps.

Finally in this section we say something about examples. One might expect
that the standard algebraic constructions for Lie algebroids
\cite{HigginsM:1990a} would have Lie bialgebroid analogues which would
provide a rich source of examples. By and large, however, this appears not
to be the case.

\begin{ex}[\cite{LiuXu:1996}]\rm                        \label{ex:LiuXu}
Let $A\to M$ be a Lie algebroid and let $\Lambda$ be an element of
$\Ga(\wedge^2A)$ with $\Lambda^\#\colon A^*\to A$ defined by
$\langle\phi,\Lambda^\#(\psi)\rangle = \Lambda(\psi, \phi)$ for
$\phi,\psi\in\Ga A^*$. Then \cite[\S2]{LiuXu:1996}, defining a bracket on
$\Ga A^*$ by
$$
[\phi,\psi] = L_{\Lambda^\#(\phi)}(\psi) - L_{\Lambda^\#(\psi)}(\phi)
               - d(\Lambda(\phi,\psi))
$$
and defining $a_* = a\smalcirc \Lambda^\#$, makes $A^*$ a Lie algebroid if
and only if $a\smalcirc[\Lambda,\Lambda]^\# = 0$ and
$L_X([\Lambda,\Lambda]) = 0$ for all $X\in\Ga A$; when these conditions
hold, it further follows that $(A,A^*)$ is a Lie bialgebroid, called
{\em exact}. It is proved in \cite{LiuXu:1996} that when $A$ is the Lie
algebroid of a Lie groupoid $G$, an exact Lie bialgebroid structure
integrates to a Poisson groupoid structure on $G$; if $\Lambda$ is the
$r$--matrix for $AG$, then the Poisson structure on $G$ is given by
$\pi = \Ri{\Lambda} - \Le{\Lambda}$. In particular, twisted Poisson groupoid
structures can be defined on the trivial Lie groupoid
$\Bar{P}\times G\times P$, where $(P,\pi_P)$
is a Poisson manifold and $G$ is a Lie group, by defining $\Lambda^{\#}\colon
T^*P\oplus(P\times{\frak g}^*)\to TP\oplus(P\times{\frak g})$ by
$\Lambda = \pi_P + X\wed\xi + r$, where $X\in{\cal X}(P)$ and
$r\colon P\to \wedge^2{\frak g}$ satisfy $L_X(\pi_P) = 0$ and
$[\xi,r] = 0$, and $[r,r]$ is pointwise ad--invariant.
\end{ex}

\begin{ex}[\cite{Kosmann-Schwarzbach:1996}]\rm
A Nijenhuis structure on a manifold $M$ is an endomorphism
$N\colon TM\to TM$
whose Nijenhuis torsion is zero. If $N$ is a Nijenhuis structure, then a
deformed Lie algebroid structure may be defined on $TM$ with $N$ as anchor
and bracket $[X,Y]_N = [NX,Y] + [X,NY] - N[X,Y]$. Denote this Lie algebroid
by $(TM)_N$.

If $M$ also has a Poisson structure $\pi$, then $N$ and $\pi$ are called
{\em compatible} if $N\pi = \pi N^*$ and a certain torsion--like expression
in $N$ and $\pi$ vanishes; $M$ with $\pi$ and $N$ is then called a
{\em Poisson--Nijenhuis structure} \cite{Kosmann-SchwarzbachM:1990}. When
this is the case, $T^*M$ with the usual cotangent Lie algebroid structure
forms a Lie bialgebroid with $(TM)_N$. Indeed compatibility of $\pi$ and
$N$ is equivalent to $T^*M$ and $(TM)_N$ forming a Lie bialgebroid.
\end{ex}

\begin{ex}[\cite{Mackenzie:SQ2}]\rm
Let $S$ be any double Lie groupoid with side groupoids $H\gpd M$ and
$V\gpd M$ and core groupoid $C\gpd M$ \cite{BrownM:1992},
\cite{Mackenzie:1992}. The horizontal groupoid structure $S\gpd V$ yields a
Lie algebroid $A_HS\to V$ which is also a groupoid over $AH$; the structure
maps of this groupoid structure are Lie algebroid morphisms, and the
resulting double structure has core Lie algebroid $AC\to M$. It follows
therefore from the duality of Pradines that $A^*_HS$ has a groupoid
structure on base $A^*C$. With respect to the dual Poisson structure,
this is a Poisson groupoid, inducing on $A^*C$ the dual Poisson structure.
Performing the same construction with the vertical structure yields a
Poisson groupoid structure on $A^*_VS\gpd A^*C$, which is dual to $A^*_HS$.

In a similar way any \LAgpd\ $(\Omega;A, V; M)$ with core Lie algebroid
$K\to M$ gives rise to a Poisson groupoid $\Omega^*\gpd K^*$.
\end{ex}

\section{Integration of Lie bialgebroids}
\label{sect:ilb}

This section is devoted to the proof of the following theorem, which is the
main result of the paper.

\begin{thm}
\label{thm:main}
Let $(AG, A^*G)$ be a Lie bialgebroid where $AG$ is the Lie algebroid of
an $\alpha$-simply connected Lie groupoid $G\gpd P$. Then there is a unique
Poisson structure on $G$ that makes $G$ into a Poisson groupoid with Lie
bialgebroid $(AG, A^*G)$.
\end{thm}

The hypotheses of the Theorem are fixed throughout the remainder of the
section.

By assumption,
\begin{equation}
\label{eq:Pi}
\Pi = \pi^\#_{AG}\smalcirc R\colon T^*A^*G\lon TAG
\end{equation}
is a morphism of Lie algebroids over $a_*$, and so
\begin{equation}
j_G\smalcirc\Pi\smalcirc s^{-1} =
               j_G\smalcirc\pi^\#_{AG}\smalcirc j'_G\colon AT^*G\lon ATG
\end{equation}
is also. Since $G$ is $\alpha$-simply connected, the cotangent groupoid
is also $\alpha$-simply connected, and so by Theorem \ref{thm:lift-morphism},
the morphism integrates uniquely to a global Lie groupoid morphism
\begin{equation}
 \matrix{&&\pi^\#&&\cr
	 &T^*G&\vlra&TG&\cr
		 &&&&\cr
			 &\vgpd&&\vgpd&\cr
				 &&&&\cr
                &A^*G&\vlra&TP&\cr
						 &&a_*&&\cr}
\end{equation}
Therefore,
\begin{equation}
\label{eq:api}
A(\pi^{\#})= j_G\smalcirc\pi^\#_{AG}\smalcirc j'_G.
\end{equation}

\begin{lem}
The map $\pi^\#$ commutes with the bundle projections,
$p_G\smalcirc \pi^\# = c_G$.
\end{lem}

\pf
Each of $j_G,\ \Pi = \pi^\#_{AG}\smalcirc R$ and $s$ is a morphism of double
vector bundles and, in particular, is a morphism of vector bundles over
$AG$. The bundle projections $AT^*G\lon AG$ and $ATG\lon AG$ are
respectively $A(c_G)$ and $A(p_G)$ so it follows that
$$
A(p_G)\smalcirc A(\pi^\#) = A(c_G).
$$
Now $p_G\smalcirc \pi^\# = c_G$ follows from the uniqueness in
Theorem \ref{thm:lift-morphism}.
\qed

\begin{lem}
\label{lem:linear}
$\pib\colon T^*G \lon TG$ is a linear map.
\end{lem}

\pf
It first has to be shown that
\begin{equation}                         \label{eq:addition}
\matrix{&&\pib \times \pib&&\cr
	&T^*G*_{G}T^*G   &\vlra&  TG*_{G}TG&\cr
		&&&&\cr
		     +&\Bigg\downarrow&&\Bigg\downarrow&+\cr
			     &&&&\cr
				     &  T^*G &\vlra&   TG&\cr
					     &&  \pib&&\cr}
\end{equation}
commutes, where both additions are the usual ones. Again, everything is a
morphism of groupoids, so we can apply $A$ and get the diagram
\begin{equation}                         \label{eq:Aaddition}
\matrix{&&A(\pib)\times A(\pib)&&\cr
   &AT^*G*_{AG}AT^*G   &\vlra&  ATG*_{AG}ATG&\cr
		&&&&\cr
           A(+)&\Bigg\downarrow&&\Bigg\downarrow&A(+)\cr
			     &&&&\cr
                 &  AT^*G &\vlra&   ATG.&\cr
                    &&  A(\pib)&&\cr}
\end{equation}
This commutes because each of $s,\ \Pi$ and $j_G$ is a morphism of vector
bundles over $AG$, and so $A(\pi^\#)$ is also. The result again follows by
the uniqueness in Theorem \ref{thm:lift-morphism}. The scalar multiplication
is handled in the same way.
\qed

\begin{lem}
$\pib$ is skew-symmetric.
\end{lem}

\pf
Since $\pi^\#$ is already known to be linear by Lemma \ref{lem:linear}, it
suffices to show that the canonical pairing $F\colon T^*G *_G TG\lon\reals$
vanishes on the graph $\calh$ of $\pi^\#$. Now $\calh$ is a subgroupoid of
$T^*G *_G TG$, and is $\alpha$-connected since $T^*G$ is, so it suffices to
show that $A(F)\colon A(\calh)\subseteq AT^*G *_{AG} ATG\lon \reals$ is zero.
But for any $\nu\in AT^*G$ we have, using Proposition \ref{pro:f-cocycle}
and Equations (\ref{eq:api}) and (\ref{eq:s}),
$$
A(F)(\nu,A(\pi^\#)(\nu))
= \langle j'_G(\nu), \pi^\#_{AG}\smalcirc j'_G(\nu)\rangle
= \pi_{AG}(j_G'(\nu),j_G'(\nu)) = 0.
$$
\qed

Combining the several lemmas above, and using
Proposition  \ref{pro:D}, we have proved
the following result.

\begin{thm}
$\pib$ defines an affine bivector field $\pi$ on $G$.
\end{thm}

Our next task is to prove that $\pi$ is indeed a Poisson tensor. As usual,
first we may introduce a bracket for functions on $G$
by\footnote{Here, as in \cite{MackenzieX:1994}, $\delta$ denotes
the usual exterior differential while the symbol  $d$ is reseverd
for the Lie algebroid differential.}
$$
\{F_{1}, F_{2}\} = \langle\pib \delta F_{1}, \delta F_{2}\rangle, \ \
\mbox{for}\ \ F_{1},  F_{2}\in C^{\infty}(G).
$$
This bracket is obviously skew-symmetric. Also, a bracket between one-forms
can be similarly introduced by
\begin{equation}
\label{eq:form-bracket}
\{\theta_1, \theta_2\} = \delta(\pi(\theta_1,\theta_2)) -
(\pib\theta_2)\per\delta\theta_1 + (\pib\theta_1)\per\delta\theta_2.
\end{equation}

\begin{lem}
\label{lem:pullback}
Let $\{\cdot , \cdot \}_{P}$ denote the Poisson bracket on $P$ induced from
the Lie bialgebroid $(AG, A^*G)$ with
$\pi^\#_P = a\smalcirc a^*_* = - a_*\smalcirc a^*$ \cite[3.6]{MackenzieX:1994}. Then
for any $f_{1}, f_{2}\in C^{\infty}(P)$,
$$
\{\alpha^* f_{1}, \alpha^* f_{2}\} = \alpha^* \{f_{1}, f_{2}\}_{P},\qquad\qquad
\{\beta^* f_{1}, \beta^* f_{2}\} = -\beta^* \{f_{1}, f_{2}\}_{P}.
$$
\end{lem}

\pf
Since $\pi^\#$ is a morphism of groupoids, we have $T\alpha \smalcirc \pib
= a_{*}\smalcirc \talp$.
Hence,
\be
\{\alpha^* f_{1}, \alpha^* f_{2}\}
&=&\langle(T\alpha \smalcirc \pib )(\alpha^* \delta f_{1}), \delta f_{2}\rangle \\
&=&\langle a_{*}\smalcirc \talp (\alpha^* \delta f_{1}), \delta f_{2}\rangle \\
&=&\langle a_{*}[-a^* \delta f_{1}],  \delta f_{2}\rangle \\
&=&\langle\pib_{P}\delta f_{1}, \delta f_{2}\rangle \\
&=&\{f_{1}, f_{2}\}_{P},
\ee
where in the third equality, we have used the fact that
$\talp (\alpha^* \delta f_{1}) = -a^* \delta f_{1}$, which follows
directly from the definition (\ref{eq:T*G}). The other identity can be
proved similarly.
\qed

Since $\pi$ is affine, $L_{\stackrel{\rightarrow}{X}}\pi$ is right invariant
for any right invariant vector field $\Ri{X}$ by Theorem \ref{thm:affine}.
The following proposition explicitly describes this invariant bivector field
(compare \cite[3.1]{Xu:1995}).

\begin{pro}
\label{pro:lie-derivative}
(i). For any  $X\in \Gamma AG$, we have
\begin{equation}
\label{eq:inf}
L_{\ssri}\pi = [\Ri{X},\pi] = -\Ri{d_{*}X},
\end{equation}
where $\Ri{d_{*}X}$ is the right invariant bivector field on $G$ corresponding
to $d_{*}X\in \seccg{2}$.

(ii). In general, for any $K\in \gm (\wedge^k AG)$, we have
$$
[\Ri{K}, \pi] = -\Ri{d_{*}K}.
$$
\end{pro}

We need a lemma to prove this proposition.

\begin{lem}
\label{lem:flow1}
Let $X\in \Gamma AG$ be any section, and let $\ell_{X}$ be the corresponding
fibrewise linear function on $A^*G$. Then the hamiltonian flow of
$\tbet^* \ell_{X}$ on $T^* G$ is $(L_{\exp{-tX}})^*$.
\end{lem}

\pf
First recall that for any vector field $X$ on any manifold $M$, if $\phi_{t}$
is a (local) flow for $X$ and $f$ is the function on $T^*M$ defined by
$f(\omega) = X\per\omega$, then the hamiltonian flow on $T^*M$ generated by
$f$ is $(\phi^{-1}_{t})^{*}$. (See Corollary 4.2.11 in \cite{AbrahamM}.)

Now, from the definition (\ref{eq:T*G}) it follows that
$(\tbet^*\ell_{X})(\omega) = \langle X,\tbet\omega\rangle  =
\langle\Ri{X},\omega\rangle$, for all $\omega\in T^* G$. Since $\Ri{X}$
has flows of the form $L_{\exp tX}$,
this completes the proof.
\qed

\nin {\bf Proof of Proposition \ref{pro:lie-derivative}}
Since $L_\ssri(\pi)$ is known to be right invariant, it suffices to
calculate $L_\ssri(\pi)(\phi_1,\phi_2)$ for $\phi_1,\phi_2\in A^*_mG$,
where we identify $\phi\in A^*_mG$ with the identity element
$\til1_\phi\in T^*_{1_m}G$ defined by $\til1_\phi(T(1)(x) + X) = \phi(X)$
for $x\in TP,\ X\in AG$. Let $\exp{tX}\in \calg (G)$ be the
family of bisections on $G$ generated by $X$, and write $g_{t} = \exp{tX}(m)$.
Now
\be
\langle L_\ssri(\pi)(1_m),\phi_1\wed\phi_2\rangle
& = & \left\langle\ddt{(L_{\exp{-tX}}\pi(g_t))},
                                          \phi_1\wed\phi_2\right\rangle\\
& = & \ddt{F(L^{*}_{\exp{-tX}}\phi_2 ,
                                    \pi^\#_{g_t}L^*_{\exp{-tX}}\phi_1)}
\ee
where $F$ is the canonical pairing. By the preceding Lemma, we have
$$
\ddt{L^{*}_{\exp{-tX}}\phi_2} = X_{\tbet^* \ell_{X}}|_{\phi_{2}}
=s(\delta\ell_X)(\phi_2),
$$
the second equality following from the definition of $s$.
Since $\pi^\#$ is a groupoid morphism, we also have
$$
\ddt{\pi^\#_{g_t}(L^*_{\exp{-tX}}\phi_1)} =
                                    A(\pi^\#)(s(\delta\ell_X)(\phi_1)).
$$
Altogether, we now have
$$
\langle L_\ssri(\pi)(1_m),\phi_1\wed\phi_2\rangle =
A(F)(s(\delta\ell_X)(\phi_2), A(\pi^\#)(s(\delta\ell_X)(\phi_1))).
$$
Using Proposition \ref{pro:f-cocycle}, this becomes
\be
\langle j_G'(s(\delta\ell_X)(\phi_2)),
         j_G^{-1}(A(\pi^\#)(s(\delta\ell_X)(\phi_1)))\rangle
& = & \langle R(\delta\ell_X(\phi_2)),\pi_{AG}^\#(R(\delta\ell_X(\phi_1)))\rangle\\
& = & \pi_{AG}(R(\delta\ell_X(\phi_1)), R(\delta\ell_X(\phi_2)))\\
& = & -d_*(X)(\phi_1,\phi_2),
\ee
the first equality following from Equations (\ref{eq:s}) and (\ref{eq:api}),
and the last equality following from \cite[6.5]{MackenzieX:1994}.

Finally, (ii) is an easy consequence of (i).
\qed

Recall that any $\omega\in\gm A^*G$ can be extended by right translation
to a linear form $\Ri{\omega}\colon T^\alpha G\to\reals$. We refer to any
1-form $\Tilde{\omega}$ on $G$ which extends $\Ri{\omega}$ as an
{\em extension} of $\omega$.

\begin{cor}
\label{cor:coin0}
Let $\Tilde{\omega}, \Tilde{ \theta}\in \Omega^{1}(G)$ be any
extensions of any two given sections $\omega, \theta\in \gm (A^*G)$,
being considered as conormal vectors. Then,
$$
\{\Tilde{\omega}, \Tilde{\theta}\}|_{P} = [\omega, \theta],
$$
where the bracket on the left hand side is defined by Equation
(\ref{eq:form-bracket}), and that on the right hand side is the Lie algebroid
bracket for $A^*G$.
\end{cor}

\pf
It follows from the definition of the Lie derivative that for any
$X\in \gm (AG)$,
$$
\Ri{X}\per \{\Tilde{\omega}, \Tilde{ \theta}\} =
(L_{\ssri}\pi )(\Tilde{\omega}, \Tilde{ \theta}) +
\pib \Tilde{\omega}\per \delta (\Ri{X}\per \Tilde{\theta} ) -
\pib \Tilde{ \theta}\per \delta (\Ri{X}\per \Tilde{\omega}).
$$
Therefore,
\be
X\per \{\Tilde{\omega}, \Tilde{ \theta}\}|_{P} & = &
(L_{\ssri}\pi )|_P (\omega, \theta)+(a_{*}\omega)  (\theta X)-(a_* \theta )
(\omega X)\\
& = & -(d_{*}X) (\omega, \theta)+(a_{*}\omega ) (\theta X)-(a_* \theta )
(\omega X)\\
& = & X\per [\omega ,\theta ],
\ee
using the definition of the Lie algebroid coboundary $d_*$. This completes the
proof.
\qed

In particular, we have

\begin{cor}
\label{cor:coin}
Let $F_{1}, F_{2}$ be any functions on $G$ which are constant
on $P$. Then $\delta\,\{F_{1}, F_{2}\}$ is conormal to $P$ and its evaluation
on $P$ is equal to $[dF_{1}, dF_{2}]$, where $d$ is the Lie algebroid
coboundary for $AG$ and the bracket is the Lie algebroid bracket on
$\gm (A^*G)$.
\end{cor}

\begin{pro}
\label{pro:schouten1}
$$
[\pi ,\pi ]|_{P}=0.
$$
\end{pro}

We need a couple of lemmas before the proof of this proposition.

\begin{lem}
\label{lem:bracket-beta}
For any $f\in C^{\infty}(P)$ and $H\in C^{\infty}(G)$,
$$
\{\beta^{*}f, H\} = -\Ri{(a_{*}^* \delta\, f)}(H).
$$
\end{lem}

\pf Commencing as in the proof of Lemma \ref{lem:pullback},
\be
\{\beta^{*}f, H\}&=&-\langle\beta^{*}\delta\, f, \pib \delta\, H\rangle \\
& = & -\langle \delta\, f, T\beta \pib \delta\, H\rangle \\
& = & -\langle \delta\, f, a_{*}\tbet \delta\, H\rangle \\
& = & -\langle a_{*}^* \delta\, f,  \tbet \delta\, H\rangle \\
& = & -\langle(\Ri{a_{*}^* \delta\, f)}, \delta\, H\rangle \\
& = & -\Ri{(a_{*}^* \delta\, f)}(H).
\ee
where, in the second last line, we used the definition (\ref{eq:T*G})
of $\tbet$. This completes the proof.
\qed

The following result summarizes the behaviour of $\pib$ on the unit space.

\begin{lem}
\label{pro:summary}
(i). For any $f\in C^{\infty}(P)$ and $m\in P$,
$\pib_{1_m}( \beta^{*}\delta\, f)$ is tangent to the $\alpha$-fiber and
equal to $-a_{*}^*\delta\, f$;

(ii). For any $\phi\in A^*G,\ \pib(\til1_\phi)$ is tangent to $P$ and
equal to $a_{*}\phi$.
\end{lem}

\nin {\bf Proof of Proposition \ref{pro:schouten1}}
It suffices to show that
\begin{equation}
\label{eq:Jacobi-P}
[\pi ,\pi]|_{P}(\delta\, F_{1}, \delta\, F_{2}, \delta\, F_{3})
= \{\{F_{1}, F_{2}\}, F_{3}\} + c.p.
\end{equation}
vanishes on $P$ for any functions $F_{1}, F_{2}, F_{3}$ defined on a
neighborhood of $P$.

Since the cotangent space $T_{1_m}^*G$ at any point $m\in P$ is spanned by
the conormal space $A^*G$ and $\beta^* T_{m}^* P$, it suffices to prove
(\ref{eq:Jacobi-P}) when the $F_{i}$ are either constant on $P$ or equal to
$\beta^* f,\ f\in C^{\infty}(P)$. We accordingly divide our proof into four
different cases.

\nin {\bf Case 1.} All $F_{i}, i=1, 2,3,$ are constant on $P$. Then
$\pib \delta\, F_{1}$ is tangent to $P$. Hence,
$\{F_{1}, F_{2}\} = \langle\pib \delta\, F_{1}, \delta\, F_{2}\rangle = 0$ on
$P$. Therefore $\{\{F_{1}, F_{2}\}, F_{3}\}$ vanishes on $P$, and so do the
other two terms.

\nin {\bf Case 2.} All $F_{i}, i=1, 2,3,$ are the pull backs of functions
on the base space $P$ by $\beta$. In this case, the conclusion follows from
Lemma \ref{lem:pullback} and the Jacobi identity of the Poisson bracket on $P$.

\nin {\bf Case 3.} Two of the $F_{i}$ are constant on $P$ and the third one
is the pull back of a function on $P$ by $\beta$. For example, assume that
$F_{1},  F_{2}$ are constant on $P$ and $F_{3} = \beta^* f_{3}$ for some
$f_{3}\in C^{\infty}(P)$. Then,
$$
\{\{F_{3}, F_{1}\}, F_{2}\} = -(a_{*}\delta\, F_{2}) \{F_{3}, F_{1}\}|_{P}
 =  (a_{*}\delta\, F_{2})[(a_{*}\delta\, F_{1} )f_{3}].
$$
Similarly, $\{\{F_{2}, F_{3}\}, F_{1}\}
= -(a_{*}\delta\, F_{1})[(a_{*}\delta\, F_{2} )f_{3}]$. On the other hand,
using Corollary \ref{cor:coin},
$$
\{\{F_{1}, F_{2}\}, F_{3}\} = a_{*}[\delta\, \{F_{1}, F_{2}\}|_{P}] (f_{3})
 = [a_{*}\delta\, F_{1}, a_{*}\delta\, F_{2}] (f_{3})
$$
The desired identity follows immediately.

\nin {\bf Case 4.} Assume that $F_{1}=\beta^* f_{1}, F_{2}=\beta^* f_{2}$
and $F_{3}$ is constant on $P$. Then,
\be
\{\{F_{2}, F_{3}\}, F_1 \} =  \{\{\beta^* f_{2}, F_{3}\}, \beta^* f_{1}\}
 = \Ri{(a^*_* \delta\, f_{1})} \{\beta^* f_{2}, F_{3}\}
 = -\Ri{(a^*_* \delta\, f_{1})}\Ri{(a^*_* \delta\, f_{2})}(F_{3}),
\ee
where we have used Lemma \ref{lem:bracket-beta} twice. Similarly,
$\{\{F_{3}, F_{1}\}, F_2 \} =
\Ri{(a^*_* \delta\, f_{2})}\Ri{(a^*_* \delta\, f_{1})}(F_{3})$.

Therefore, using \cite[3.7]{MackenzieX:1994},
$$
\{\{F_{2}, F_{3}\}, F_1 \}+\{\{F_{3}, F_{1}\}, F_2 \}
 = [\Ri{(a^*_* \delta\, f_{2})},\Ri{(a^*_* \delta\, f_{1})}](F_{3})
 = \Ri{[a^*_* \delta\, \{f_{2}, f_{1}\} ]}(F_{3}).
$$
On the other hand,
$$
\{\{F_{1}, F_{2}\}, F_3 \} =
\{\{ \beta^* f_{1}, \beta^* f_{2}\}, F_{3}\}
 = -\{ \beta^* \{f_{1},f_{2}\}, F_{3}\}
 = \Ri{[a^*_* \delta\, \{f_{1}, f_{2}\} ]}(F_{3}).
$$
Thus the RHS of Equation (\ref{eq:Jacobi-P}) vanishes on $P$.
This completes the proof.
\qed

\nin {\bf Proof of Theorem \ref{thm:main}}  Let $D=[\pi, \pi]$. Then $D$
is an affine bivector field according to Theorem \ref{thm:affine}. For
any $X\in \gm (AG)$, $L_{\ssri}D = 2[L_{\ssri}\pi , \pi] =
-2 [\Ri{(d_{*}X)}, \pi] = 2\Ri{(d_{*}^2 X)} = 0$ according to Proposition
\ref{pro:lie-derivative}. Thus, $dD=0$. Since $D|_P = 0$ by Proposition
\ref{pro:schouten1}, it follows that $D=0$ according to
Theorem \ref{thm:dD0}. That is, $\pi$ is indeed a Poisson tensor.
Therefore $G$ is a
Poisson groupoid according to Proposition 8.1 in \cite{MackenzieX:1994}.
Finally, Corollary \ref{cor:coin0} implies that the associated Lie bialgebroid
is exactly isomorphic to $(AG,A^*G)$. The uniqueness of the Poisson structure
is quite evident again by \cite{MackenzieX:1994}.
\qed

Note that the integrability of exact Lie bialgebroids proved in
\cite{LiuXu:1996} and described above in Example~\ref{ex:LiuXu} is more
precise than Theorem \ref{thm:main} and does not require any
connectedness hypotheses.

\section{Symplectic groupoids}
\label{sect:sg}

The notion of symplectic groupoid plays an important role in Poisson
geometry. There already exists an extensive body of work on this
subject---see, for example, \cite{CDW}, \cite{DazordW:1991}, \cite{KarasevM}.
One of the main reasons for studying Poisson groupoids and Lie bialgebroids
is to understand better the relation between symplectic groupoids and Poisson
groups by putting both of them into this more general framework. A particular
problem is to understand how the symplectic structure of a symplectic groupoid
arises. Symplectic groupoids are, of course, Poisson groupoids whose Poisson
structure is symplectic. Thus one main difference between Poisson groups and
general Poisson groupoids is that Poisson groups are never symplectic but the
latter may be symplectic. It turns out that there is a close relation between
symplectic groupoids and their base Poisson manifolds. In fact, they are in
one-one correspondence in a rough local sense \cite{CDW},
\cite{Karasev:1989}, \cite{Karasev:1987}.

Let $P$ be a Poisson manifold with Poisson tensor $\pi_{P}$. It is
well-known that the cotangent bundle $T^{*}P\lon P$ carries a natural Lie
algebroid structure.
Given $f\in \cinf(P)$, denote the Hamiltonian vector field corresponding to
$f$ by $X_f$. Then the anchor $\pi^\#\colon T^*P\lon TP$ is determined
by $\pi^\#(f\delta\, g) = fX_g$. Given $\omega, \theta\in\Omega^{1}(P)$, and
writing $X_{\omega}=\pi^\#\omega$ and $X_{\theta} = \pi^\#\theta$, the Lie
algebroid bracket is
\begin{equation}
\label{eq:cotangent}
\{\omega, \theta\}=L_{X_{\omega}}\theta - L_{X_{\theta}}\omega -
\delta\,(\pi (\omega, \theta)).
\end{equation}
On the other hand, the tangent bundle $TP$ of $P$ has the trivial Lie
algebroid structure given by the usual bracket of vector fields.

\begin{pro}
{\em (i)}. Let $P$ be a Poisson manifold, and $T^*P$ the cotangent Lie
algebroid described above. Then $(T^* P, TP)$ is a Lie bialgebroid.

{\em (ii)}. Conversely, any Lie algebroid structure on $T^*P$ which is
compatible with the trivial tangent bundle Lie algebroid $TP$ in the sense
that they become a Lie bialgebroid arises in this way.
\end{pro}

\pf  See Example 3.3 in \cite{MackenzieX:1994} for (i).

Let $(A, A^{*})$ be a Lie bialgebroid such that $A^*$ is isomorphic to $TP$ as
a Lie algebroid. It is easy to see that this isomorphism must be realized by
the anchor $a_{*}: A^*\lon TP$. According to \cite[3.6]{MackenzieX:1994},
there is a Poisson structure on $P$ induced from the Lie bialgebroid
$(A, A^{*})$. Hence $T^*P $ has a Lie algebroid structure as defined by
Equation (\ref{eq:cotangent}), and $(T^*P, TP)$ becomes a Lie bialgebroid.
Furthermore, the anchor  $a_{*}: A^*\lon TP$ is in fact a Lie bialgebroid
morphism. Hence, its dual $a_{*}^* :T^*P\lon A$ is a Lie algebroid
morphism. Since $a_{*}$ is an isomorphism, so is $a_{*}^*$. In other
words, $A$ is isomorphic to $T^*P$ arising from a Poisson structure on $P$.
\qed

Our main theorem of the section is the following.

\begin{thm}
\label{thm:symplectic}
Suppose that $P$ is a Poisson manifold and the corresponding cotangent Lie
algebroid $T^*P$ integrates to an $\alpha$-simply connected groupoid $G$.
Then $G$ admits a natural symplectic structure which makes it into a
symplectic groupoid.
\end{thm}

This result follows from our integration Theorem \ref{thm:main} and the
following:

\begin{thm}
Let $\poidd{G}{P}$ be a Poisson groupoid, with Lie bialgebroid $(AG, A^*G)$.
Let $\pi$ denote the Poisson tensor on $G$. Then the following are equivalent:
\begin{enumerate}
\item $(AG, A^*G)$ is isomorphic to the canonical Lie bialgebroid
$(T^*P , TP)$ associated to the base Poisson manifold $P$;
\item $\pi$ is non-degenerate along $P$;
\item $G$ is a symplectic groupoid.
\end{enumerate}
\end{thm}

\pf $\mbox{(i)}\Longleftrightarrow \mbox{(ii)}.$
Fix any  point  $m\in P$. The cotangent space $T_{1_m}^*G$ is spanned by the
conormal space $A^*_mP$ and by the covectors of the form $\beta^* \delta\, f$
for $f\in C^{\infty}(P)$, so any $\omega\in T_{1_m}^*G$ can be written as
$\phi + \beta^* \delta\, f$ for some $\phi \in A^*_mP$ and $f\in C^{\infty}(P)$.
Now suppose that $\pib (\phi + \beta^* \delta\, f) = 0$. Then, it follows from
Proposition \ref{pro:summary} that $a_{*}\phi - a_*^* \delta\, f = 0$, which
implies immediately that $a_{*}\phi = 0$ and $a_*^* \delta\, f = 0$, since
$a_{*}\phi$ is tangent to $P$ while $a_*^* \delta\, f$ is tangent to the
$\beta$-fiber. From this, we conclude that $\pib_{1_m}$ is nondegenerate if
and only if both
$a_{*}: A^*_{m}G\lon T_{m}P$ and $a_{*}^* :T^*_{m}P\lon A_{m}G$ are
one-to-one. This, however, is exactly equivalent to saying that
$a_{*}: A^*_{m}G\lon T_{m}P$ is one-to-one and onto. Hence, that $\pi$ is
non-degenerate along $P$
is equivalent to saying that $a_{*}: A^*G\lon TP$ is a vector bundle
isomorphism, and this is also equivalent to it being a Lie bialgebroid
isomorphism since it is already a Lie bialgebroid morphism.
This completes our proof of the first part.

$\mbox{(ii)}\Longleftrightarrow \mbox{(iii)}.$
We only need to prove one direction, namely, that if $\pi$ is nondegenerate
along $P$, then it is nondegenerate everywhere. Let $X\in \gm_{c} (AG)$ be
any compactly supported section which is closed in the sense that
$d_{*}X = 0$. Then, according to Corollary 3.6 in \cite{Xu:1995},
$\exp{X}$ is a coisotropic bisection, and therefore the
corresponding right translation  $R_{\exp{X}}: G\lon G$ is a Poisson
diffeomorphism. Hence $\pi$ is nondegenerate along the bisection $\exp{X}$.
According to the first part of the proof, the Lie bialgebroid is in fact
isomorphic to $(T^*P, TP)$. Therefore there exist abundant closed sections
for $AG = T^*P$ in the sense that through any point of $A$ there exists
a closed section.  Hence, through any point of $G$, there exists a bisection
which is a product of those of the form $\exp{X}$ for some closed section.
In other words, there always exists a coisotropic bisection through any point
of $G$. Therefore, $\pi$ is nondegenerate everywhere.
\qed

{\bf Remark:} Theorem \ref{thm:symplectic} fails in general if $G$ is not
assumed to be $\alpha$-simply connected. An example due to Weinstein is
given in \cite[5.2]{Dazord:1994}.

According to Pradines \cite{Pradines:4}, any Lie algebroid can be integrated
to a local groupoid, which can of course be assumed to be $\alpha$-simply
connected. As an immediate consequence of Theorem \ref{thm:symplectic}, we
therefore obtain the following theorem of Karasev \cite{Karasev:1987} and
Weinstein \cite{Weinstein:1987} on the existence of local symplectic
groupoids for arbitrary Poisson manifolds.

\begin{cor}
There always exists a local symplectic groupoid over any Poisson manifold.
In particular, any Poisson manifold admits a symplectic realization.
\end{cor}

Although the groupoid structure is only locally defined, it includes
the entire base manifold and so gives a global symplectic realization.

\section{Appendix: Lifting of Lie algebroid morphisms}
\label{sect:app}

This section is devoted to the proof of the following theorem.

\begin{thm}
\label{thm:lift-morphism}
Let $\poidd{G}{M}$ and $\poidd{H}{N}$ be Lie groupoids. Suppose that $G$ is
$\alpha $-simply connected, and that $\phi: AG\lon AH,\ \phi_{0}: M\lon N$
is a Lie algebroid morphism. Then $\phi$ integrates uniquely to a Lie
groupoid morphism $\Phi: G\lon H$.
\end{thm}

\pf
Let $\Delta \subseteq M\times N$ be the graph of $f$, and consider
$G*H = (\beta_{G}\times \beta_H)^{-1}(\Delta)$. Also, let $B$ be the
graph of $\phi$. Then $B$ is a subalgebroid with base $\Delta$ of the
product Lie algebroid $AG\times AH$. Identifying $B$ with the tangent
subbundle of $\alpha$-fibers along the base $\Delta$ and translating it
by right translations, one obtains a distribution $\cald$ on $G*H$, which
is integrable since $B$ is a subalgebroid.

Given any $m\in M$, let $\call_{m}$ be the leaf of $\cald$ through the
point $(1_m, 1_{f(m)})$. Consider the projection
$p_G\colon \call_{m}\lon G$. It is clear, since $\cald$ is defined using
right translations, that $p_G( \call_{m})\subseteq G_{m}$, where
$G_{m}=\alpha_G^{-1}(m)\subseteq G$. We will show that
$p_G\colon \call_{m}\lon G_{m}$ is a covering map. For this, it suffices
to prove that $p_G$ possesses the path lifting property.

It is simple to see that the section space $\gm (AG)$ can be naturally
identified with the section space $\gm (B)$, which in turn can be identified
with the space of right invariant vector fields on $G*H$. More explicitly,
given any $X\in \gm (AG)$, its corresponding right invariant vector field
on $G*H$, denoted by $\Ri{X}_{B}$, is given by
$$
\Ri{X}_{B}(g, h) = (T(R_{g})X_{n}, T(R_{h})\phi (X_{n})),
                                       \ \ \forall (g, h)\in G*H,
$$
where $n=\beta_G(g)$ and $X_{n}\in AG|_{n}$. At the same time,
$X$ defines a right invariant vector field on $G$, which is denoted by
$\Ri{X}$. Let $\tau_{t}$ and $\rho_{t}$ denote flows of $\Ri{X}_{B}$
and $\Ri{X}$, respectively. Using Lemma~\ref{lem:KS} below, it is simple
to see that $\tau_{t}(g, h)$ is defined for all $(g, h)\in G*H$ and
$t\in \reals$, whenever $\rho_{t}(g)$ is defined for all $g$ and $t$.
Since $p_{1*}\Ri{X}_{B}=\Ri{X}$, it also follows that
$$
p_1\smalcirc \tau_{t}(g, h)=\rho_{t}(g).
$$

Suppose that $\sigma_{t}$ is any path in $G_{m}$ starting from $1_m$ and
which is a product of paths generated by right invariant vector fields
$\Ri{X_{i}}$, for compactly supported sections $X_{i}\in \gm (AG)$. Let
$\Tilde{\sigma}_{t}$ be the corresponding path in $G*H$, the product
of the flows generated by the right invariant vector fields $\Ri{X_{i}}_{B}$
in the same order. Clearly, $\Tilde{\sigma}_{t}$ lies in $\call_{m}$, and is
a lift of $\sigma_{t}$. Since any path in $G_{m}$ can be approximated
by such paths $\sigma_{t}$, this proves the path lifting property
for $p_G$.

Since $G_{m}$ is simply connected by assumption, it follows that
$p_G\colon \call_{m}\lon G_{m}$ is a diffeomorphism. Therefore it defines
a smooth map $\Phi_{m}\colon G_{m}\lon H_{f(m)}$, where
$H_{f(m)}=\alpha^{-1}_H(f(m))\subseteq H$. Since $\call_{m}\subseteq G*H$,
$\Phi_{m}$ commutes with $\beta_G$ and $\beta_H$. In fact, $\Phi_{m}$ is
characterized by the following properties:

\begin{enumerate}
\item $\Phi_m(1_m)= 1_{f(m)}$;
\item $\beta_H\smalcirc \Phi_{m} = \beta_G$;
\item $T_g(\Phi_m)(Y) = T(R_{\Phi_m(g)})[\phi(T(R_{g^{-1}}(Y)))], \ \
\mbox{for all}\ g\in G_{m}\ \mbox{ and } Y\in T_{g}G_{m}$.
\end{enumerate}

By varying $m$ in $M$, we obtain a family of smooth maps $\Phi_{m}:
G_{m}\lon H_{f(m)}$, and hence a global map $\Phi$ from $G$ to $H$.

Fix any $y\in G$ and suppose that $m=\beta_G(y)$, $n=\alpha_G(y)$.
Consider the map $\Psi : G_{m}\lon H_{f(m)}$ given by
$$
\Psi (x)=\Phi_{n}(xy)\Phi^{-1}_{n}(y), \ \ \forall x\in G_{m}.
$$
It is simple to see that $\Psi$ satisfies the properties (i)--(iii) above.
Therefore, $\Psi=\Phi_{m}$, which is equivalent to
$$
\Phi_{n}(xy)=\Phi_{m}(x)\Phi_{n}(y).
$$
In other words, $\Phi$ is a groupoid morphism.

It remains to prove that $\Phi$, as a global map from $G$ to $H$, is smooth.
First, it is easy to see that, by construction, it is smooth along the
identity space $M$. In fact, the graph of $\Phi$, in a neighborhood of any
$m\in M$, is the flow box in $G*H$ obtained from flowing the graph of $f$
along the distribution $\cald$, which is transversal to it; thus we obtain
a smooth submanifold.

Now take any point $x_0\in G$. As in Lemma~\ref{lem:enough}, there exist
$X_1,\ldots,X_n\in\Ga_c(AG)$ such that $\exp X_1\ldots\exp X_n$ has the
value $x_0$ at $\alpha x_0$. For notational convenience in what follows,
assume that $n = 1$, so that $X$ has $\exp X(\alpha x_0) = x_0$.

\begin{lem}
\label{lem:lambda}
The map $\lambda\colon M\to H$ defined by
$\lambda (m) = \Phi(\exp{X}(m)),\ m\in M,$ is smooth.
\end{lem}

\pf
Let $\gamma$ be the graph map of $f$ regarded as $M\to G*H$; thus
$\gamma(m) = (1_m,1_{f(m)})$. It is easy to see, from the definition, that
$\lambda$ is the composition $p_{2}\smalcirc \tau_1\smalcirc\gamma$, where
$\tau_{1}$ is the time-$1$ map of the flow $\tau_{t}$ introduced early in the
proof, and $p_{2}: G*H\lon H$ is the projection. Hence $\lambda$ is smooth.
\qed

Since $\Phi$ is a groupoid morphism, for $x$ in a neighborhood of $x_{0}$,
we have
$$
\Phi(x) = \Phi (\exp{X}(m(x)) \cdot \exp{X}(m(x))^{-1} \cdot x)
 =  \Phi(\exp{X}(m(x)))\Phi(\exp{X}(m(x))^{-1}\cdot x)
$$
where $m(x)=\Ad_{\exp^{-1}X}(\beta_G x)$; note that $\Ad_{\exp X}$
is the time-1 map for $a(X)\in {\cal X}(M)$.

According to Lemma \ref{lem:lambda}, $\Phi(\exp{X}(m(x)))$ is a smooth map
from $G$ to $H$. On the other hand,
$\exp X(m(x))^{-1} = \exp -X(\beta_G x)$, and so
$\Phi(\exp{X}(m(x))^{-1}\cdot x) =
(\Phi \smalcirc L_{\exp -X})(x)$ is also smooth in a neighborhood of $x_{0}$.
Since the groupoid multiplication is smooth, it follows that $\Phi$ is
smooth in a neighborhood of $x_{0}$. This concludes the proof of
Theorem~\ref{thm:lift-morphism}.
\qed

Finally we recall the result of Kumpera and Spencer used above.

\begin{lem}[\cite{KumperaS}]
\label{lem:KS}
Let $G$ be any Lie groupoid over base $M$, and let $AG$ be its Lie algebroid
with anchor $a$. For $X\in\Ga AG$, let $\Ri{X}$ be the right invariant vector
field corresponding to $X$, and recall that $a(X)= \beta_{*}\Ri{X}$ is its
projected vector field on $M$. Then $\Ri{X}$ is complete if and only if
$a(X)$ is complete. In fact, $\Tilde{\sigma}_{t}(x)$ is defined whenever
$\sigma_{t}(\beta (x))$ is defined, where $\Tilde{\sigma}_{t}$ and
$\sigma_{t}$ are the flows generated by $\Ri{X}$ and $a(X)$, respectively.
\end{lem}

Theorem \ref{thm:lift-morphism} was announced by Pradines \cite{Pradines:1};
a method of proof (based on the local integrability of Lie subalgebroids)
was briefly indicated by Almeida and Kumpera \cite{AlmeidaK:1981}.
When both Lie algebroids are transitive and over the same base, with the
map $f$ an identity, the proof can be reduced to standard results of
connection theory: see \cite[III\S7]{Mackenzie:LGLADG}, which also gives
further references for this case. Providing the target Lie algebroid is
transitive, the case of a general base map may be reduced to that of an
identity by using a pullback \cite{HigginsM:1990a}. When the target is a
Lie algebra (so that the base map is constant), the result reduces to the
integration of Maurer--Cartan forms: see \cite[\S5]{Xu:1995}; this method
of proof can be extended to handle any transitive target. For the case of
a compact base (and $f$ an identity), a proof has been given by Mokri
\cite{Mokri:1996}. See also \cite{BrownMu:1995} for the closely related
construction of the monodromy groupoid of an $\alpha$--connected groupoid.

\newcommand{\noopsort}[1]{} \newcommand{\singleletter}[1]{#1}


\begin{thebibliography}{10}

\bibitem{AbrahamM}
R.~Abraham and J.~Marsden.
\newblock {\em Foundations of Mechanics}.
\newblock Addison-Wesley, second edition, 1985.

\bibitem{AlbertD:1991ssr}
C.~Albert and P.~Dazord.
\newblock Th{\'e}orie des groupo{\"\i}des symplectiques: {C}hapitre {II},
  {G}roupo{\"\i}des symplectiques.
\newblock In {\em Publications du D{\'e}partement de Math{\'e}matiques de
  l'Universit{\'e} {C}laude {B}ernard, Lyon {I}}, nouvelle s{\'e}rie, pages
  27--99, 1990.

\bibitem{AlmeidaK:1981}
R.~Almeida and A.~Kumpera.
\newblock Structure produit dans la cat{\'e}gorie des alg{\`e}broids de {L}ie.
\newblock {\em An. Acad. Brasil Ci{\^e}nc}, 53:247--250, 1981.

\bibitem{BrownM:1992}
R.~Brown and K.~C.~H. Mackenzie.
\newblock Determination of a double {L}ie groupoid by its core diagram.
\newblock {\em J.~Pure Appl. Algebra}, 80(3):237--272, 1992.

\bibitem{BrownMu:1995}
R.~Brown and O.~Mucuk.
\newblock The monodromy groupoid of a {L}ie groupoid.
\newblock {\em Cahiers Topologie G{\'e}om. Diff{\'e}rentielle
  Cat{\'e}goriques}, 36(4):345--369, 1995.

\bibitem{CDW}
A.~Coste, P.~Dazord, and A.~Weinstein.
\newblock Groupo{\"\i}des symplectiques.
\newblock In {\em Publications du D{\'e}partement de Math{\'e}matiques de
  l'Universit{\'e} de Lyon, {I}}, number 2/A-1987, pages 1--65, 1987.

\bibitem{Dazord:1994}
P.~Dazord.
\newblock Lie groups and {L}ie algebras in infinite dimension: {A} new
  approach.
\newblock In Y.~Maeda, H.~Omori, and A.~Weinstein, editors, {\em Symplectic
  geometry and quantization}, number 179 in Contemporary Mathematics, pages
  17--44. Amer. Math. Soc., July 1993.

\bibitem{DazordW:1991}
P.~Dazord and A.~Weinstein, editors.
\newblock {\em Symplectic geometry, groupoids and integrable systems},
  S{\'e}minaire {S}ud {R}hodanien de {G}{\'e}om{\'e}trie (1989).
  Springer-Verlag, MSRI Publications, 20, 1991.

\bibitem{Drinfeld:1983}
V.~G. Drinfel'd.
\newblock Hamiltonian structures on {L}ie groups, {L}ie bialgebras and the
  geometric meaning of the classical {Y}ang-{B}axter equation.
\newblock {\em Soviet. Math. Dokl.}, 27:68--71, 1983.

\bibitem{HigginsM:1990a}
P.~J. Higgins and K.~C.~H. Mackenzie.
\newblock Algebraic constructions in the category of {L}ie algebroids.
\newblock {\em J.~Algebra}, 129:194--230, 1990.

\bibitem{Karasev:1987}
M.~V. Karasev.
\newblock Analogues of the objects of {L}ie group theory for nonlinear
  {P}oisson brackets.
\newblock {\em Math. USSR-Izv.}, 28:497--527, 1987.

\bibitem{Karasev:1989}
M.~V. Karasev.
\newblock The {M}aslov quantization conditions in higher cohomology and analogs
  of notions developed in {L}ie theory for canonical fibre bundles of
  symplectic manifolds.~{I}, {II}.
\newblock {\em Selecta Math. Soviet.}, 8:213--234, 235--258, 1989.
\newblock Preprint, Moscov. Inst. Electron Mashinostroeniya, 1981, deposited at
  VINITI, 1982.

\bibitem{KarasevM}
M.~V. Karasev and V.~P. Maslov.
\newblock {\em Nonlinear {P}oisson brackets: {G}eometry and {Q}uantization},
  volume 119 of {\em Translations of Mathematical Monographs}.
\newblock American Mathematical Society, Providence, R.I., 1993.

\bibitem{Kosmann-Schwarzbach:1995}
Y.~Kosmann-Schwarzbach.
\newblock Exact {G}erstenhaber algebras and {L}ie bialgebroids. {G}eometric and
  algebraic structures in differential equations.
\newblock {\em Acta Appl. Math.}, 41:153--165, 1995.

\bibitem{Kosmann-Schwarzbach:1996}
Y.~Kosmann-Schwarzbach.
\newblock The {L}ie bialgebroid of a {P}oisson--{N}ijenhuis manifold.
\newblock {\em Lett. Math. Phys.}, 38:421--428, 1996.

\bibitem{Kosmann-SchwarzbachM:1990}
Y.~Kosmann-Schwarzbach and F.~Magri.
\newblock Poisson-{N}ijenhuis structures.
\newblock {\em Ann. Inst. H.~Poincar{\'e} Phys. Th{\'e}or.}, 53:35--81, 1990.

\bibitem{KumperaS}
A.~Kumpera and D.~C. Spencer.
\newblock {\em Lie equations. {V}olume I: {G}eneral theory}.
\newblock Princeton University Press, 1972.

\bibitem{Lie}
S.~Lie.
\newblock {\em Theorie der {T}ransformationsgruppen. {Z}weiter {A}bschnitt,
  unter {M}itwirkung von {P}rof. {D}r. {F}riedrich {E}ngel}.
\newblock Teubner, Leipzig, 1890.

\bibitem{LiuWX:Dirac}
Zhang-Ju Liu, Alan Weinstein, and Ping Xu.
\newblock Dirac structures and {P}oisson homogeneous spaces.
\newblock Preprint, University of California, Berkeley, 1996.

\bibitem{LiuWX:1997}
Zhang-Ju Liu, Alan Weinstein, and Ping Xu.
\newblock Manin triples for {L}ie bialgebroids.
\newblock {\em J. Differential Geom.}, 45:547--574, 1997.

\bibitem{LiuXu:1996}
Zhang-Ju Liu and Ping Xu.
\newblock Exact {L}ie bialgebroids and {P}oisson groupoids.
\newblock {\em Geom. Funct. Anal.}, 6:138--145, 1996.

\bibitem{LuW:1990}
Jiang-Hua Lu and A.~Weinstein.
\newblock Poisson {L}ie groups, dressing transformations, and {B}ruhat
  decompositions.
\newblock {\em J. Differential Geom.}, 31:501--526, 1990.

\bibitem{Mackenzie:LGLADG}
K.~Mackenzie.
\newblock {\em Lie groupoids and {L}ie algebroids in differential geometry}.
\newblock London Mathematical Society Lecture Note Series, no.~124. Cambridge
  University Press, {\noopsort{1985}}1987.

\bibitem{Mackenzie:1992}
K.~C.~H. Mackenzie.
\newblock Double {L}ie algebroids and second-order geometry, {I}.
\newblock {\em Adv. Math.}, 94(2):180--239, 1992.

\bibitem{Mackenzie:SQ2}
K.~C.~H. Mackenzie.
\newblock Double {L}ie algebroids and second-order geometry, {II}.
\newblock Preprint, \noopsort{1997}1996.

\bibitem{MackenzieX:1994}
K.~C.~H. Mackenzie and Ping Xu.
\newblock Lie bialgebroids and {P}oisson groupoids.
\newblock {\em Duke Math.~J.}, 73(2):415--452, 1994.

\bibitem{MackenzieX:1998}
K.~C.~H. Mackenzie and Ping Xu.
\newblock Classical lifting processes and multiplicative vector fields.
\newblock {\em Quarterly J. Math.}, 1998.
\newblock To appear.

\bibitem{Mokri:1996}
T.~Mokri.
\newblock On {L}ie algebroid actions and morphisms.
\newblock {\em Cahiers Topologie G{\'e}om. Diff{\'e}rentielle
  Cat{\'e}goriques}, 37:315--331, 1996.

\bibitem{Pradines:1}
J.~Pradines.
\newblock Th{\'e}orie de {L}ie pour les groupo{\"\i}des diff{\'e}rentiables.
  {R}elations entre propri{\'e}t{\'e}s locales et globales.
\newblock {\em C.~R. Acad. Sci. Paris, S{\'e}rie~A}, 263:907--910, 1966.

\bibitem{Pradines:4}
J.~Pradines.
\newblock Troisi{\`e}me th{\'e}or{\`e}me de {L}ie pour les groupo{\"\i}des
  diff{\'e}rentiables.
\newblock {\em C.~R. Acad. Sci. Paris, S{\'e}rie~A}, 267:21--23, 1968.

\bibitem{Pradines:1988}
J.~Pradines.
\newblock Remarque sur le groupo{\"\i}de cotangent de {W}einstein-{D}azord.
\newblock {\em C.~R. Acad. Sci. Paris S{\'e}r.~I Math.}, 306:557--560, 1988.

\bibitem{Weinstein:1983}
A.~Weinstein.
\newblock The local structure of {P}oisson manifolds.
\newblock {\em J. Differential Geom.}, 18:523--557, 1983.
\newblock Errata and Addenda, same journal, 22:255, 1985.

\bibitem{Weinstein:1987}
A.~Weinstein.
\newblock Symplectic groupoids and {P}oisson manifolds.
\newblock {\em Bull. Amer. Math. Soc. (N.S.)}, 16:101--104, 1987.

\bibitem{Weinstein:1988}
A.~Weinstein.
\newblock Coisotropic calculus and {P}oisson groupoids.
\newblock {\em J. Math. Soc. Japan}, 40:705--727, 1988.

\bibitem{Weinstein:1990}
A.~Weinstein.
\newblock Affine {P}oisson structures.
\newblock {\em Internat. J. Math.}, 1:343--360, 1990.

\bibitem{Xu:1995}
Ping Xu.
\newblock On {P}oisson groupoids.
\newblock {\em Internat. J.~Math.}, 6(1):101--124, 1995.

\end{thebibliography}
\end{document}